\documentclass{ar-1col}

\usepackage{graphicx}
\usepackage{amsmath}
\usepackage{amsthm}
\usepackage{amsfonts}
\usepackage{amssymb}
\usepackage{graphicx}
\usepackage{enumerate}
\usepackage{caption}
\usepackage{subfigure}
\usepackage{wrapfig}
\usepackage{color}
\usepackage{algpseudocode}
\usepackage{algorithm}
 \usepackage[active]{srcltx}
\usepackage[unicode]{hyperref}
\usepackage[numbers]{natbib} 
\usepackage{hyperref}
\usepackage{soul}
\usepackage{mathtools}
\usepackage{xcolor}

\newtheorem{definition}{Definition}
\newtheorem{proposition}{Proposition}
\newtheorem{theorem}{Theorem}

\newcommand*{\ani}{{}}

\newcommand*{\aaa}{{}}

\jname{Xxxx. Xxx. Xxx. Xxx.}
\jvol{AA}
\jyear{YYYY}
\doi{10.1146/((please add article doi))}

\begin{document}

\markboth{Majumdar, Hall, Ahmadi}{Scalability in Semidefinite Programming}

\title{A Survey of Recent Scalability Improvements for Semidefinite Programming
with Applications in Machine Learning, Control, and Robotics}

\author{Anirudha Majumdar,$^1$ Georgina Hall,$^2$ and Amir Ali Ahmadi$^3$
\affil{$^1$Department of Mechanical and Aerospace Engineering, Princeton University, Princeton, USA, 08544; email: ani.majumdar@princeton.edu}
\affil{$^2$Decision Sciences, INSEAD, Fontainebleau, France, 77300; email: georgina.hall@insead.edu}
\affil{$^3$Department of Operations Research and Financial Engineering, Princeton University, Princeton, USA, 08544; email: a$\_$a$\_$a$\_$@princeton.edu}}

\begin{abstract}
\vspace{-5pt}
Historically, scalability has been a major challenge to the successful application of semidefinite programming in fields such as machine learning, control, and robotics. In this paper, we survey recent approaches for addressing this challenge including (i) approaches for exploiting structure (e.g., sparsity and symmetry) in a problem, (ii) approaches that produce low-rank approximate solutions to semidefinite programs, (iii) more scalable algorithms that rely on augmented Lagrangian techniques and the alternating direction method of multipliers, and (iv) approaches that trade off scalability with conservatism (e.g., by approximating semidefinite programs with linear and second-order cone programs). For each class of approaches we provide a high-level exposition, an entry-point to the corresponding literature, and examples drawn from machine learning, control, or robotics. We also present a list of software packages that implement many of the techniques discussed in the paper. Our hope is that this paper will serve as a gateway to the rich and exciting literature on scalable semidefinite programming for both theorists and practitioners.  
\end{abstract}

\maketitle

\tableofcontents

\section{INTRODUCTION}
\label{sec:intro}

A \emph{semidefinite program} (SDP) is an optimization problem of the form
\begin{equation}\label{eq:SDP}
\begin{aligned}
&\min_{X \in S_n} &&\mbox{Tr}(CX)\\
&\text{s.t. } &&\mbox{Tr}(A_iX)=b_i,~i=1,\ldots,m\\
& &&X \succeq 0,
\end{aligned}
\end{equation}
where $S_n$ denotes the set of $n\times n$ real symmetric matrices, $\mbox{Tr}(.)$ stands for the trace of a matrix, $C, A_1,\ldots, A_m\in S_n$ and $b_1,\ldots,b_m\in \mathbb{R}$ are input data to the problem, and the notation $X \succeq 0$ denotes that the matrix $X$ is constrained to be \emph{positive semidefinite} (psd); i.e., to have nonnegative eigenvalues. The \emph{dual} to Problem \ref{eq:SDP} takes the form
\begin{equation}\label{eq:dual SDP}
\begin{aligned}
&\max_{y \in \mathbb{R}^m} &&\sum_{i=1}^m y_i b_i\\
&\text{s.t. } &&C - \sum_{i=1}^m y_i A_i \succeq 0 
\end{aligned}
\end{equation}
\aaa{and will also feature in this paper.}

SDPs have numerous fundamental applications throughout applied and computational mathematics. In combinatorics and discrete optimization for instance, some of the most powerful relaxations for prominent problems of the field rely on semidefinite programming. For example, the current best approximation algorithm for the maximum-cut problem~\cite{goemans1995improved}, or the only known polynomial-time algorithm for finding maximum stable sets in perfect graphs~\cite{lovasz1979shannon} make fundamental use of semidefinite programming.
In machine learning, one is often faced with the task of learning a signal with a known structure from data (e.g., a low-rank matrix, a sparse eigenvector, a monotone polynomial, etc.). The design of an appropriate convex optimization problem which induces such structure in a signal often leads to a semidefinite programming problem; see, e.g.,~\cite{BenMaryamPabloRank},~\cite{dAspremont07},~\cite[Chap. 8]{GH_thesis}.
In polynomial optimization, SDP-based hierarchies based on the notion of sum of squares polynomials (see Section~\ref{subsec:sos} below) are among the most promising approaches for obtaining global solutions in the absence of any convexity assumptions~\cite{lasserre_moment,sdprelax}. Finally, in control and robotics, various notions of safety, performance, stability, and robustness of dynamical systems can be certified by searching for Lyapunov functions that satisfy appropriate nonnegativity constraints over the state space. When the candidate functions are parameterized as polynomials or other semialgebraic functions, this search can be automated by semidefinite programming via the sum of squares approach; see, again, Section~\ref{subsec:sos} below and~\cite{PhD:Parrilo}. Overall, semidefinite programming has proven to be one of the most expressive classes of convex optimization problems, subsuming many other important classes such as linear programming, convex quadratic programming, and second-order cone programming~\cite{VaB:96}.




In addition to its expressiveness, another contributing factor to the popularity of semidefinite programming lies in its strong theoretical and computational properties, which to certain extent parallel that of linear programming~\cite{VaB:96}. For instance, strong duality between Problem \ref{eq:SDP} and Problem \ref{eq:dual SDP} holds under mild assumptions.
Moreover, many algorithms for linear programs, such as the ellipsoid method and interior point methods, can and have been adapted to SDPs. 
While these algorithms solve SDPs to arbitrary accuracy in polynomial time, in practice, they suffer from one serious impediment: \emph{scalability}. When the problems considered have large $n$ or $m$ (as defined in Problem \ref{eq:SDP}), solving time and memory requirements tend to explode, making SDPs prohibitive to \aaa{work with}. Indeed, most solvers rely on interior point methods as their algorithm of choice for SDPs and interior point methods produce iterates that are (as their name suggests) in the interior of the feasible set. Computing these large dense matrix iterates and storing them lead to the issues described above. Devising methods to make SDPs more scalable is consequently a very important and active area of research. 

\subsection{The goal and outline of this survey}

In this paper, we review recent literature on scalability improvements for semidefinite programming with a particular focus on methods that have proved useful in machine learning, control theory, and robotics. We consider these representative application areas since they are timely and because many problems that arise in them have straightforward semidefinite programming-based relaxations or reformulations, but involve very large problem instances. This paper is a follow-up to a workshop at the 2016 Conference on Decision and Control on the same topic and  titled \emph{``Solving Large SDPs in Control, Machine Learning, and Robotics''}. Eleven lectures were given at this workshop whose content can be found at \url{http://aaa.princeton.edu/largesdps}. 

The outline of this paper is as follows. In Section \ref{subsec:sos}, we provide a brief introduction to sum of squares proofs of nonnegativity of polynomials. This background material is relevant to many of the applications of semidefinite programming we discuss in this paper (e.g., in Section \ref{sec:structure pop} and Section~\ref{subsec:dsos}). Section \ref{sec:structure} presents approaches that enhance scalability by exploiting structure (e.g., symmetry or sparsity) in a problem. Section \ref{sec:low.rank} presents methods for generating low-rank solutions to SDPs in order to reduce computation time and storage requirements. Section \ref{sec:first order} reviews more scalable alternatives to interior-point methods for solving SDPs based on augmented Lagrangian techniques and the alternating direction method of multipliers. Section \ref{sec:scalability_conservatism} presents approaches that trade off scalability with conservatism of solutions (e.g., by approximating SDPs with linear or second-order cone programs). In Section \ref{sec:software}, we present a list of software packages that implement many of the techniques described in this paper. Section \ref{sec:conclusions} concludes the paper and highlights promising directions for future work.

\subsubsection{Relationship between approaches discussed in the paper} \ani{The approaches for improving scalability we discuss in Sections \ref{sec:structure}--\ref{sec:scalability_conservatism} are largely complementary to each other. These sections in the paper can thus be read independently of each other. We provide the following rough guide to the different sections (which may be used, e.g., by a practitioner interested in a specific application where scalability is a challenge). The interested reader may also use the guide below for exploring the software packages listed in Section \ref{sec:software}\aaa{, which} are also organized according to sections in the paper. 
\begin{itemize}
\item If the SDP of interest (corresponding to Problem \ref{eq:SDP}) has special \emph{structure} (e.g., symmetry, sparsity, or degeneracy) $\rightarrow$ Section \ref{sec:structure}.
\item If either (i) the SDP has \emph{low-rank} (optimal) solutions, or (ii) one desires good-quality low-rank feasible points to the SDP  $\rightarrow$  Section \ref{sec:low.rank}.
\item If \emph{approximately feasible} solutions with good objective values are of interest (i.e., slight violation of the constraints can be tolerated) $\rightarrow$ Section \ref{sec:first order}.
\item If \emph{convservative feasible} solutions to the SDP  (i.e., points that satisfy all constraints, but are potentially suboptimal) are valuable $\rightarrow$ Section \ref{sec:scalability_conservatism}.
\end{itemize}
In principle, the approaches presented in the different sections may be combined with each other to further enhance scalability (e.g., in order to tackle a problem exhibiting sparsity for which slightly suboptimal feasible solutions are valuable). We highlight such combinations of approaches where possible in Sections \ref{sec:structure}--\ref{sec:scalability_conservatism}. }

\subsection{Related work not covered by this survey}
The number of algorithms that have been proposed as alternatives to interior point methods for semidefinite programming is large and rapidly growing. While we have attempted to present some of the major themes, our survey is certainly not exhaustive.  Some of the interesting contributions that this paper does not cover due to space limitations include the recent first-order methods developed by Renegar~\cite{renegar2019accelerated}, the multiplicative weights update algorithm of Arora, Hazan, and Kale~\cite{arora2005fast},  the storage-optimal algorithm of Ding et al.~\cite{ding2019optimal}, the iterative projection methods of Henrion and Malick~\cite{henrion2011projection},  the conditional gradient-based augmented Lagrangian framework of Yurtsever, Fercoq, and Cevher~\cite{yurtsever2019conditional}, the regularization methods of Nie and Wang for SDP relaxations of polynomial optimization~\cite{nie2012regularization}, the accelerated first-order method of Bertsimas, Freund, and Sun for the sum of squares relaxation of unconstrained polynomial optimization~\cite{bertsimas2013accelerated}, and the spectral algorithms of Hopkins et al. for certain sum of squares problems arising in machine learning~\cite{hopkins2016fast}.

We also remark that there has been a relatively large recent literature at the intersection of optimization and machine learning where SDP relaxations of certain nonconvex problems are avoided altogether and instead the nonconvex problem is tackled directly with local descent algorithms. These algorithms can often scale to larger problems compared to the SDP approach (at least when the SDPs are solved by interior point methods). Moreover, under certain technical caveats and statistical assumptions on problem data, it is sometimes possible to prove similar guarantees for these algorithms as those of the SDP relaxation. We do not cover this literature in this paper but point the interested reader to~\cite{nonconvex_list} for a list of references.


\subsection{A brief review of sum of squares proofs of nonnegativity}\label{subsec:sos}
Sum of squares constraints on multivariate polynomials are responsible for a substantial fraction of recent applications of semidefinite programming.
For the convenience of the reader, we briefly review the basics of this concept and the context in which it arises. This will be relevant e.g. for a better understanding of Section \ref{sec:structure pop} and Section~\ref{subsec:dsos}.

Recall that a \emph{closed basic semialgebraic set} is a subset of the Euclidean space of the form
\begin{equation}\nonumber
\Omega\mathrel{\mathop:}=\{x\in\mathbb{R}^n|\ g_i(x)\geq 0, i=1,\ldots,m  \},
\end{equation}
where $g_1,\ldots,g_m$ are polynomial functions. It is relatively well known by now that many fundamental problems of computational mathematics can be cast as optimization problems where the decision variables are the coefficients of a multivariate polynomial $p$, the objective function is a linear function of these coefficients, and the constraints are twofold: (i) affine constraints in the coefficients of $p$, and (ii) the constraint that 
\begin{equation}\label{eq:p>=0.on.Omega}
p(x)\geq 0\  \forall x\in\Omega,
\end{equation}
where $\Omega$ is a given closed basic semialgebraic set. For example, in a \emph{polynomial optimization problem}---i.e., the problem of minimizing a polynomial function $f$ on $\Omega$---the optimal value is equal to the largest constant $\gamma$ for which $p(x)\mathrel{\mathop:}=f(x)-\gamma$ satisfies Constraint \ref{eq:p>=0.on.Omega}. See~\cite{PabloGregRekha_BOOK},\cite{GH_Business},\cite[Sect. 1.1]{dsos_siaga} for numerous other applications. Unfortunately, imposing Constraint \ref{eq:p>=0.on.Omega} (or even asking a decision version of it for a fixed polynomial $p$) is NP-hard already when $p$ is a quartic polynomial and $\Omega=\mathbb{R}^n$, or when $p$ is quadratic and $\Omega$ is a polytope.

An idea pioneered to a large extent by Lasserre~\cite{lasserre_moment} and Parrilo~\cite{sdprelax} has been to write algebraic sufficient conditions for Constraint \ref{eq:p>=0.on.Omega} based on the concept of sum of squares polynomials. We say that a polynomial $h$ is a \emph{sum of squares} (sos) if it can be written as $h=\sum_i q_i^2$ for some polynomials $q_i$. Observe that if we succeed in finding sos polynomials $\sigma_0,\sigma_1,\ldots, \sigma_m$ such that the polynomial identity
\begin{equation}\label{eq:Putinar.proof}
p(x)=\sigma_0(x)+\sum_{i=1}^{m}\sigma_i(x)g_i(x)
\end{equation}
holds, then, clearly, Constraint \ref{eq:p>=0.on.Omega} must be satisfied. Conversely, a celebrated result in algebraic geometry~\cite{putinar1993positive} states that if $g_1,\ldots,g_m$ satisfy the so-called ``Archimedean property'' (a condition slightly stronger than compactness of the set $\Omega$), then positivity of $p$ on $\Omega$ guarantees existence of sos polynomials $\sigma_0,\sigma_1,\ldots, \sigma_m$ that satisfy the algebraic identity in Equation \ref{eq:Putinar.proof}.


The computational appeal of the sum of squares approach stems from the fact that the search for sos polynomials $\sigma_0,\sigma_1,\ldots, \sigma_m$ of a given degree that verify the polynomial identity in Equation \ref{eq:Putinar.proof} can be automated via semidefinite programming. This is true even when some coefficients of the polynomial $p$ are left as decision variables. This claim is a straightforward consequence of the following well-known fact (see, e.g.,~\cite{PhD:Parrilo}): A polynomial $h$ of degree $2d$ is a sum of squares if and only if there exists a symmetric matrix $Q$ which is positive semidefinite and verifies the identity $$h(x)=z(x)^TQz(x),$$ where $z(x)$ here denotes the vector of all monomials in $x$ of degree less than or equal to $d$. Note that the size of the semidefinite constraint in an SDP that imposes a sum of squares constraint on an $n$-variate polynomial of degree $2d$ is $ {n+d \choose d} \times  {n+d \choose d}\approx n^d\times n^d.$ This is the reason why \emph{sum of squares optimization problems}---i.e., SDPs arising from sum of squares constraints on polynomials---often run into scalability limitations quickly.



\section{ENHANCING SCALABILITY BY EXPLOITING STRUCTURE}
\label{sec:structure}

A major direction for improving the scalability of semidefinite programming has been towards the development of methods for exploiting problem-specific structure. Existing literature in this area has focused on two kinds of structure that appear frequently in applications in control theory, robotics, and machine learning: (i) sparsity, and (ii) symmetry.



\subsection{Exploiting sparsity}
\label{sec:sparsity}

Many problems in control theory give rise to SDPs that exhibit structure in the form of \emph{chordal sparsity}. We first give a brief introduction to the general notion of chordal sparsity and then highlight applications in control theory that exploit this structure to achieve significant gains in scalability. We begin by introducing a few relevant graph-theoretic notions. 


\begin{definition}[Cycles and chords]
A \emph{cycle} of length $k \geq 3$ in a graph $\mathcal{G}$ is a sequence of distinct vertices $(v_1, v_2, \dots, v_k)$ such that $(v_k, v_1)$ is an edge and $(v_i, v_{i+1})$ is an edge for $i = 1, \dots, k-1$. A \emph{chord} (associated with a cycle) is an edge that is not part of the cycle but connects two vertices of the cycle.  
\end{definition}

\begin{definition}[Chordal graph]
An undirected graph $\mathcal{G}$ is \emph{chordal} if every cycle of length $k \geq 4$ has a chord. 
\end{definition}

\begin{definition}[Maximal clique]
A \emph{clique} of a graph $\mathcal{G}$ is a subset $\mathcal{C}$ of vertices such that for any distinct vertices $i,j \in \mathcal{C}$, $(i,j)$ is an edge in $\mathcal{G}$. A clique is \emph{maximal} if it is not a subset of another clique.  
\end{definition}

The graph theoretic notions introduced above can be used to capture sparsity patterns of matrices. Let $\mathcal{E}$ (resp. $\mathcal{V}$) denote the set of edges (resp. vertices) of an undirected graph $\mathcal{G}$. 
Let $S_n$ denote the set of symmetric $n \times n$ matrices as before, and $S_n(\mathcal{E})$ denote the set of matrices with a sparsity pattern defined by $\mathcal{G}$ as follows:
\begin{equation}
S_n(\mathcal{E}) \coloneqq \{X \in S_n \ | \  X_{ij} = 0 \ \text{if} \ (i,j) \notin \mathcal{E} \ \mbox{for} \  i \neq j \}.
\end{equation}

The following proposition is the primary result that allows one to exploit chordal sparsity (i.e., sparsity induced by a chordal graph) in order to improve scalability of a semidefinite program. 

\begin{proposition}[see \cite{Agler88}]
\label{prop:agler}
Let $\mathcal{G}$ be a chordal graph with vertices $\mathcal{V}$, edges $\mathcal{E}$, and a set of maximal cliques $\Gamma = \{\mathcal{C}_1, \dots, \mathcal{C}_p\}$. Then $X \in S_n(\mathcal{E})$ is positive semidefinite if and only if there exists a set of matrices $X_k \in S(\mathcal{C}_k) \coloneqq  \{X \in S_n \ | \ X_{ij} = 0 \ \text{if} \ (i,j) \notin \mathcal{C}_k \times \mathcal{C}_k \}$ for $k = 1,\dots,p$, such that $X_k \succeq 0$ and $X = \sum_{k=1}^{p} X_k$. 
\end{proposition}

The above statement allows one to equivalently express a large semidefinite constraint as a set of smaller semidefinite constraints (and additional equality constraints). This can lead to a significant increase in computational efficiency since the computational cost of solving an SDP is primarily determined by the size of the largest semidefinite constraint. We note that for a chordal graph, the problem of listing maximal cliques is amenable to efficient (i.e., polynomial time) algorithms (see \cite{Zheng17} and references therein).

Proposition \ref{prop:agler} (and a related dual result due to Grone et al. \cite{Grone84}) has been leveraged by several researchers in optimization to improve the efficiency of SDPs with chordal sparsity \cite{Fukuda01, Kim11, Fujisawa09}. These results have also been exploited by researchers in the context of control theory. In \cite{Mason14}, the authors apply these results to the problem of finding Lyapunov functions for sparse linear dynamical systems. In particular, by parameterizing Lyapunov functions with a chordal structure (i.e., a structure that allows one to apply Proposition \ref{prop:agler}), the resulting SDP finds a Lyapunov function significantly faster (up to a factor of $\sim 80$ faster in terms of running time for instances where the largest maximal clique size is less than 15). 

In \cite{Zheng16, Zheng17}, the authors extend Proposition \ref{prop:agler} to the case of matrices with block-chordal sparsity by demonstrating that such matrices can be decomposed into an equivalent set of smaller block-positive semidefinite matrices (again, with additional equality constraints). This result allows for applications to networked systems since it reasons about the sparsity of interconnections between subsystems. The authors apply these results to the problem of designing structured feedback controllers to stabilize large-scale networked systems.

In the work highlighted above, Proposition \ref{prop:agler} is used to decompose a large semidefinite constraint into smaller ones with additional equality constraints. The resulting SDP can then be solved using standard techniques (e.g., standard interior point methods). One challenge with this approach (as noted in \cite{Zheng19}) is that the additional equality constraints can sometimes nullify the computational benefits of dealing with smaller semidefinite constraints. A set of strategies that seek to overcome this issue involve exploiting chordal sparsity \emph{directly in the solver}, e.g., directly in an interior point method \cite{Fukuda01, Burer03, Andersen10} or in a first-order method \cite{Sun14, Kalbat15, Madani15, Sun15, Zheng19}. In the context of control theory, \cite{Andersen14} leverages sparse solvers \cite{Andersen10, Fukuda01, Benson06} based on this idea to tackle robust stability analysis for large-scale sparsely interconnected systems. 

\subsubsection{Exploiting sparsity in polynomial optimization problems}
\label{sec:structure pop}

The existence of structure in the form of sparsity can also be exploited for polynomial optimization problems (cf. Section \ref{subsec:sos}). In \cite{Weisser18}, the authors propose the Sparse-BSOS hierarchy (which is based on the Bounded-SOS (BSOS) hierarchy for polynomial optimization problems presented in \cite{Lasserre17}). The Sparse-BSOS hierarchy leverages the observation that for many polynomial optimization problems that arise in practice, the polynomials that define the constraints of the problem exhibit sparsity in their coefficients. This observation is used to split the variables in the problem into smaller blocks of variables such that (i) each monomial of the polynomial objective function only consists of variables from one of the blocks, and (ii) each polynomial that defines the constraints also only depends on variables in only one of the blocks. Under the assumption that the sparsity in the objective and constraints satisfy the Running Intersection Property (RIP) \cite{Weisser18}, the Sparse-BSOS approach provides a hierarchy of SDPs whose optimal values converge to the global optimum of the original polynomial optimization problem. In addition, for the so-called class of ``SOS-convex polynomial optimization problems" (see \cite{Weisser18} for a definition) that satisfy RIP, the hierarchy converges at the very first step. In contrast to the (standard) BSOS hierarchy (which also provides these guarantees), Sparse-BSOS results in semidefinite constraints of smaller block size. Numerical experiments in \cite{Weisser18} demonstrate that exploiting sparsity in the problem can result in the ability to solve large-scale polynomial optimization problems that are beyond the reach of the standard (i.e., non-sparse) BSOS approach.  

In \cite{Mangelson18}, the authors utilize the Sparse-BSOS hierarchy to tackle the problem of \emph{simultaneous localization and mapping (SLAM)} in robotics. SLAM \cite{Thrun05} refers to the problem of simultaneously estimating the location (or trajectory of locations) of a robot and a model of the robot's environment. Two important versions of the SLAM problem are (i) Landmark SLAM, and (ii) Pose-graph SLAM. The Landmark SLAM problem involves simultaneously estimating the position and orientation of the robot and the location of landmarks in the environment. The Pose-graph SLAM problem is a more restricted version that only requires estimating the position and orientation of the robot at each time-step. Traditionally, these problems have been posed as maximum likelihood estimation (MLE) problems and tackled using general nonlinear optimization methods (which do not come with guarantees on finding globally optimal solutions). In \cite{Mangelson18}, the authors formulate the planar pose-graph and landmark SLAM problems as polynomial optimization problems and demonstrate that they are SOS-convex. This allows \cite{Mangelson18} to apply the Sparse-BSOS hierarchy with the guarantee that the hierarchy converges to the global optimum at the very first step. The approach in \cite{Mangelson18} contrasts with most prior works on pose-graph/landmark SLAM, which apply general-purpose nonlinear optimization algorithms to these problems. The numerical results in \cite{Mangelson18} compare the Sparse-BSOS approach to SLAM with nonlinear optimization techniques based on the Levenberg-Marquardt algorithm and demonstrate that the latter often yield suboptimal solutions while the former finds the globally optimal solution. We also note that the approach in \cite{Mangelson18} contrasts with the SE-Sync approach \cite{Rosen19} we mention in Section \ref{subsec:BM} since \cite{Mangelson18} relaxes the requirement of limited measurement noise imposed by \cite{Rosen19}. 

\subsection{Exploiting symmetry}
\label{sec:symmetry}

Beyond sparsity, another general form of structure that arises in applications of semidefinite programming to control theory, robotics, and machine learning is \emph{symmetry}. Mathematically, symmetry refers to a transformation of the variables in a problem that does not change the underlying problem. As a concrete example, consider a multi-agent robotic system composed of a large number of identical agents. The agents in the system can be interchanged while leaving the overall capabilities of the multi-agent system invariant.  

\emph{Symmetry reduction} \cite{Gatermann04, Vallentin09, deKlerk10, Permenter17} is a powerful set of approaches for exploiting symmetries in SDPs. We provide a brief overview of symmetry reduction here and refer the reader to \cite{Permenter17} for a thorough exposition. Let $S_n^+$ denote the cone of $n \times n$ symmetric positive semidefinite matrices. Given an SDP of the form of Problem \ref{eq:SDP}, symmetry reduction performs the following two steps:
\begin{enumerate}
\item Find an appropriate (see \cite{Permenter17}) affine subspace $\mathcal{S} \subset S_n$ containing feasible solutions to the problem;
\item Express $\mathcal{S} \cap S_n^+$ as a linear transformation of a ``simpler" cone $\mathcal{C} \subset S_m$, where $m < n$:
\begin{equation}
\mathcal{S} \cap S_n^+ = \{\Phi Z \ | \ Z \in \mathcal{C} \}.
\end{equation}
\end{enumerate} 
Symmetry reduction techniques work by taking $\mathcal{S}$ equal to the \emph{fixed point subspace} of the group action corresponding to the underlying symmetry group (see \cite{Permenter17, Gatermann04} for details). 

Once steps 1 and 2 above have been accomplished, one can solve the following problem instead of the original SDP:
\begin{equation}\label{eq:SDP_symmetry}
\begin{aligned}
&\min_{Z \in S_m} &&\mbox{Tr}(C \Phi Z)\\
&\text{s.t. } &&\mbox{Tr}(A_i \Phi Z)=b_i,~i=1,\ldots,m\\
& &&Z \in \mathcal{C}.
\end{aligned}
\end{equation}
The key advantage here is that Problem \ref{eq:SDP_symmetry} is formulated over a lower-dimensional space $S_m$, and can thus be more efficient to solve.

Symmetry reduction techniques have been fruitfully applied in many different applications of semidefinite programming to control theory \cite{Henrion05, Cogill08, Arcak16}. As an example, \cite{Cogill08} exploits symmetry in linear dynamical systems; more precisely, it exploits invariances of the state space realizations of linear systems under unitary coordinate transformations to obtain significant computational gains for control synthesis problems. In \cite[Chapter 7]{Arcak16}, symmetry reduction techniques are used to reduce the computational burden of certifying stability of interconnected subsystems. Specifically, \cite{Arcak16} exploits permutation symmetries in a network of dynamical systems (i.e., invariances to exchanging subsystems) in order to reduce the complexity of the resulting SDPs. 

\subsection{Facial reduction}
\label{sec:facial reduction}

In certain cases, there may be additional structure beyond sparsity and symmetry that one can exploit in order to reduce the size of SDPs arising in practice. As an example, consider the search for a Lyapunov function $V: \mathbb{R}^n \rightarrow \mathbb{R}$ that proves stability of the equilibrium $x_0 \in \mathbb{R}^n$ of a polynomial dynamical system $\dot{x} = f(x)$. In this case, one would like to search for functions $V$ that satisfy the following ``degeneracy" conditions: $V(x_0) = 0, \dot{V}(x_0) = 0$. This imposes structure on the resulting optimization problem by restricting the space of possible functions $V$. One powerful set of approaches for leveraging this kind of ``degeneracy" structure is referred to as \emph{facial reduction} \cite{Permenter17, Drusvyatskiy17, Borwein81, Pataki13, kungurtsev2018two}. The general procedure behind facial reduction is identical to that of symmetry reduction (Steps 1 and 2 in Section \ref{sec:symmetry}). The primary difference is that facial reduction techniques identify the subspace $\mathcal{S}$ in a different way than symmetry reduction techniques. In particular, $\mathcal{S}$ is found by exploiting a certain kind of geometric degeneracy condition (see \cite{Permenter17} for details).

Facial reduction techniques have been used to improve scalability of SDPs arising in robotics, control theory, and ML applications \cite{Permenter17, Waki16, Krislock10}. We highlight \cite{Permenter17}, which presents methods for \emph{automating} the process of performing facial (and symmetry) reduction. As an example, \cite{Permenter17} considers the problem of analyzing the stability of a \emph{rimless wheel}, which is a simple hybrid dynamical walking robot model. The facial reduction techniques presented in \cite{Permenter17} reduce computation times by roughly a factor of $60$ for this problem. Similar gains in terms of scalability are also obtained for many other problems. Another recent facial reduction algorithm which stands out in terms of its simplicity of implementation is ``\emph{Sieve-SDP}''~\cite{sieve_sdp}. This technique can detect redundancies in the constraints of an SDP (or sometimes infeasibility) by checking positive definiteness of some sub-blocks of the SDP data matrices. In~\cite{sieve_sdp}, the authors demonstrate the performance of this preprocessing strategy on several benchmark examples, including many arising from SDP relaxations of polynomial optimization problems.

\section{ENHANCING SCALABILITY BY PRODUCING LOW-RANK SOLUTIONS}
\label{sec:low.rank}

In this section, we focus on methods for producing low-rank feasible solutions to SDP \ref{eq:SDP} at relatively large scales. One may think of the techniques we present here as leveraging a very specific kind of ``structure" (i.e., low-rank structure), similar to the approaches considered in Section \ref{sec:structure}. However, the technical approaches for leveraging low-rank structure are quite different in nature from the approaches presented in Section \ref{sec:structure}. Moreover, the literature on low-rank methods is quite vast. We thus treat these methods distinctly from the ones in Section \ref{sec:structure}. At a high level, one utilizes such methods in two different settings, which we describe now. \\

\noindent \textbf{Case 1: There is a low-rank optimal solution to Problem \ref{eq:SDP}.} In this setting, we can establish a priori that among the optimal solutions to Problem \ref{eq:SDP}, there must be one of low rank. An important case where this is known to hold is when Problem \ref{eq:SDP} does not have too many constraints. More specifically (see~\cite{barvinok1995problems,pataki1998rank,lemon2016low}), if Problem \ref{eq:SDP} has $m$ constraints and an optimal solution, then it also has a rank-$r$ optimal solution with 
\begin{align} \label{eq:low.rank.exist}
\frac{r(r+1)}{2} \leq m.
\end{align} Note that we need $m < \frac{n(n+1)}{2}$ for this property to be useful.
	
In this setting, one may not particularly care about obtaining a low-rank (optimal) solution to the problem. Rather, searching for a low-rank solution (whose existence is guaranteed) can lead to significant gains in computation time and storage. To see why, recall that a $n \times n$ positive semidefinite matrix $X$ has rank $r \leq n$ if and only if it can be written as $X=UU^T$, where $U \in \mathbb{R}^{n \times r}$. By searching for a low-rank solution, one can hope to work in the lower-dimensional space $\mathbb{R}^{n \times r}$ corresponding to $U$, rather than the higher-dimensional space $S_n$ corresponding to $X$. Doing so would lead to much cheaper algorithms as storage space needs would drop from $O(n^2)$ to $O(nr)$ and, e.g., the important operation of matrix-vector multiplication would require $O(nr)$ flops rather than $O(n^2)$. \\ 

\noindent \textbf{Case 2: Good-quality low-rank feasible solutions to Problem \ref{eq:SDP} are of interest.} In this setting, one considers SDPs whose optimal solutions may all have high rank. However, the goal is to find low-rank feasible solutions with good objective values. By making this compromise, one can again expect to work in the lower-dimensional space of low-rank feasible solutions and make gains in terms of computation time and storage as explained previously.

Problems such as these occur in a variety of applications, most notably machine learning and combinatorial optimization. Consider for instance the problem of low-rank matrix completion~\cite{candes2010matrix} from machine learning, which is a cornerstone problem in user recommendation systems (see, e.g., the Netflix competition \cite{bennett2007netflix}). In this problem, we have access to a matrix $Z \in \mathbb{R}^{m \times n}$, which is partially complete, i.e., some of its entries are specified but others are unknown. The goal is to fill in the unknown entries. Without assuming some structure on $Z$, any completion of $Z$ would be a valid guess and so the problem is ill-defined. This can be remedied by assuming that $Z$ should have low rank, which is a natural assumption in the context of user recommender systems. Indeed, if one assumes that the matrix $Z$ reflects the ratings a user $i$ ($i=1,\ldots,n$) gives a product $j$ ($j=1,\ldots,m$), then completing the matrix $Z$ would involve infering how a user would have rated a product (s)he has not yet rated from existing ratings, provided by both him/herself and other users. Under the assumption that a consumer base can be segmented into $r$ categories (with $r\ll n$), it would then make sense to constrain the matrix $Z$ to have rank less than or equal $r$, which would reflect this segmentation. The matrix completion problem would then be written as 
\begin{equation}\label{eq:matrix.completion}
\begin{aligned}
&\min_{X \in \mathbb{R}^{m \times n}} &&\sum_{i,j} (X_{ij}-Z_{ij})^2\\
&\text{s.t. } &&\mbox{rank}(X) \leq r,
\end{aligned}
\end{equation}
where $Z$ is the partially complete matrix corresponding to past data. The indices $i,j$ in the objective function range over known entries of the matrix $Z$. Problem \eqref{eq:matrix.completion} can be relaxed to an SDP by replacing the constraint on the rank of $X$ by a constraint on its \emph{nuclear norm} \cite{BenMaryamPabloRank}, which is a semidefinite-representable surrogate of the rank. The problem then becomes:
\begin{equation}\label{eq:matrix.completion.relaxation}
\begin{aligned}
&\min_{X \in \mathbb{R}^{m \times n}, W_1 \in S_m, W_2 \in S_n} \sum_{i,j} (X_{ij}-Z_{ij})^2\\
&\text{s.t. } \mbox{Tr}\begin{bmatrix} W_1 & X\\ X^T & W_2 \end{bmatrix} \leq r, \begin{bmatrix} W_1 & X\\ X^T & W_2 \end{bmatrix}  \succeq 0.
\end{aligned}
\end{equation}
The goal here is to obtain a feasible solution to the problem such that the matrix 
$$\begin{bmatrix} W_1 & X\\ X^T & W_2 \end{bmatrix}$$
has low rank and the objective value is small. This makes matrix completion a typical application of the algorithms we present next. Other machine learning problems which involve relaxing a rank-constrained SDP to an SDP via the use of the nuclear norm include certain formulations of clustering and maximum variance unfolding \cite{kulis2007fast}, and sparse principal component analysis \cite{dAspremont07}.

The rest of this section is devoted to two of the most well-known approaches for computing low-rank feasible (and possibly optimal) solutions to SDP \ref{eq:SDP}. The first is the Burer-Monteiro method and its variants (Section \ref{subsec:BM}). The second includes methods that use the Frank-Wolfe algorithm as their computational backbone (Section \ref{subsec:FW}). 

\subsection{Burer-Monteiro \aaa{based} methods} \label{subsec:BM}
The key idea of the Burer-Monteiro algorithm \cite{burer2003nonlinear} is to factor the decision variable $X$ in Problem \ref{eq:SDP} as $VV^T$ where $V \in \mathbb{R}^{n \times r}$. Here, $r$ can either be chosen such that a solution of rank $r$ is guaranteed to exist (e.g., such that Inequality \ref{eq:low.rank.exist} holds) (Case 1 in the previous segment) or such that it corresponds to the rank which we would like our feasible solution to have (Case 2). Typically, the algorithm comes with some theoretical guarantees for the first case and is simply a heuristic in the second. Using this factorization, we rewrite Problem \ref{eq:SDP} as:
\begin{equation}\label{eq:BM}
\begin{aligned}
&\min_{V \in \mathbb{R}^{n \times r}} &&\mbox{Tr}(CVV^T)\\
&\text{s.t. } &&\mbox{Tr}(A_iVV^T)=b_i, i=1,\ldots,m.
\end{aligned}
\end{equation}
This reformulation has both pros and cons: on the upside, the problem to consider has smaller dimension than Problem \ref{eq:SDP} and has no conic constraint to contend with; on the downside, the problem has become nonconvex. Hence, local optimization methods may not always recover the global minimum of Problem \ref{eq:SDP} even if there is one of low rank. Research in this area has consequently sought to find algorithms to solve Problem \ref{eq:BM} that (i) provably recover global low-rank solutions to the problem if they exist; (ii) do so reasonably fast in terms of both theoretical and empirical convergence rates; and (iii) work on a large class of inputs $A_i, b_i,$ and $C$. We review existing algorithms through the lens of these criteria.

\subsubsection{Augmented Lagrangian approaches}
The original method recommended by Burer and Monteiro \cite{burer2003nonlinear} to solve Problem \ref{eq:BM} is a method based on an augmented Lagrangian procedure (see Section \ref{subsec:Aug.Lag.Method} for another appearance of this general approach). The idea is to rewrite Problem \ref{eq:BM} as an unconstrained optimization problem by adding terms to the objective that penalize infeasible points. This is done by constructing an augmented Lagrangian function, given here by
$$L_{y,\sigma}(V)=\mbox{Tr}(CVV^T)-\sum_{i=1}^m y_i (\mbox{Tr}(A_iVV^T)-b_i)+\frac{\sigma}{2} \sum_{i=1}^m (\mbox{Tr}(A_iVV^T)-b_i)^2.$$
If $y$ and $\sigma >0$ are fixed appropriately, then minimizing $L_{y,\sigma}(\cdot)$ will yield an optimal solution to Problem \ref{eq:BM}. To obtain an appropriate pair $(y,\sigma)$, the following iterative procedure is followed: fix $(y,\sigma)$ and minimize $L_{y,\sigma}(\cdot)$ to obtain $V$, then construct a new $(y,\sigma)$ from $V$, and repeat; see \cite{burer2003nonlinear} for the exact details. Note that this procedure involves the minimization of $L_{y,\sigma}$ with respect to $V$: choosing an appropriate algorithm for this subproblem is also an important step in the method. The authors in \cite{burer2003nonlinear} suggest using a limited-memory BFGS (LM-BFGS) algorithm \cite{liu1989limited} to do so. The LM-BFGS algorithm is a quasi-Newton algorithm for unconstrained optimization problems, where the objective is differentiable. Unlike Newton methods, it does not require computing the Hessian of the objective function, constructing instead an approximate representation of it from past gradients. Given that gradients of $L_{y,\sigma}$ can be computed quite quickly under some conditions (see \cite{burer2003nonlinear} for more details), this is a reasonably efficient algorithm to use. One can replace the LM-BFGS subroutine with a subgradient approach if Problem \ref{eq:BM} has additional constraints that are non-differentiable. Furthermore, the augmented Lagrangian approach has been tailored to directly handle inequality constraints of the type $\mbox{Tr}(A_iX) \leq b_i$ in Problem \ref{eq:SDP}; see \cite{kulis2007fast}. Though there are no theoretical guarantees regarding convergence to global solutions, the authors of \cite{burer2003nonlinear} experimentally observe that the algorithm tends to return global minima of SDP \ref{eq:SDP}, doing so in a much smaller amount of time than interior point or bundle methods. The follow-up paper \cite{burer2005local} to \cite{burer2003nonlinear} slightly modifies the augmented Lagrangian function and shows convergence of the algorithm to a global minimum of SDP \ref{eq:SDP} (given that one exists of rank $\leq r$) if, among other things, each minimization of the modified $L_{y,\sigma}(\cdot)$ gives rise to a local minimum, a condition that is hard to test in practice.

\subsubsection{Riemannian optimization approaches}
Riemannian optimization concerns itself with the optimization of a smooth objective function over a smooth Riemannian manifold. We remind the reader that a \emph{Riemannian} manifold $\mathcal{M}$ is a manifold that can be linearized at each point $x \in \mathcal{M}$ by a tangent space $T_x \mathcal{M}$, which is equipped with a Riemannian metric. We say that a Riemannian manifold is \emph{smooth} if the Riemannian metric varies smoothly as the point $x$ varies; see e.g. \cite[Appendix A]{boumal2011rtrmc} for a short introduction to such concepts.  In the Riemannian setting, one can define notions of both gradients and Hessians for manifolds \cite[Appendix A]{boumal2011rtrmc}. Using these notions, one can extend the concepts of first-order and second-order necessary conditions of optimality to optimization over manifolds. These are almost identical to the ones for unconstrained optimization over $\mathbb{R}^n$, except that they involve the analogs of gradient and Hessians for manifolds. Similarly, one can generalize many classical algorithms for nonlinear unconstrained optimization (such as, e.g., trust-region and gradient-descent methods) to the manifold setting. We refer the reader to \cite{absil2007trust,absil2009optimization} for more information on the specifics of these algorithms which we will not cover here.

Though it may not be immediately apparent why, Riemannian optimization methods have become a popular tool for tackling problems of the form \ref{eq:BM} \cite{journee2010low,boumal2016non}. Indeed, the set
$$\mathcal{M}=\{V \in \mathbb{R}^{n \times r}~|~\mbox{Tr}(A_iVV^T)=b_i, i=1,\ldots,m\}$$
is a smooth Riemannian manifold under certain conditions on $A_i$ (see \cite[Assumption 1]{boumal2018deterministic}), and the objective function $V \mapsto \mbox{Tr}(CVV^T)$ is smooth. Hence, Problem \ref{eq:BM} is exactly a Riemannian optimization problem. It can be shown that if $m < \frac{r(r+1)}{2}$, the conditions on $\mathcal{M}$ that make it smooth hold, and $\mathcal{M}$ is compact, then for almost all $C \in S_n$, any second-order critical point of Problem \ref{eq:BM} is globally optimal to Problem \ref{eq:SDP} \cite[Theorem 2]{boumal2016non}. (Note that an optimal solution of rank $\leq k$ is guaranteed to exist here from the assumptions and that we use the manifold definition of a second-order critical point.) If an algorithm which returns second-order critical points can be exhibited, then, in effect, we will have found a way of solving Problem \ref{eq:BM}. It turns out that Riemannian trust-region methods return such a point---under some conditions---regardless of initialization \cite{boumal2018global}. In fact, it is shown in \cite{boumal2018global} that an $\epsilon$-accurate second-order critical point (i.e., a point $V \in \mathcal{M}$ with $||\mbox{grad}(\mbox{Tr}(CVV^T))|| \leq \epsilon$ and $\mbox{Hessian}(\mbox{Tr}(CVV^T)) \succeq -\epsilon I$ ) can be obtained in $O(1/\epsilon^2)$ iterations. 
Other Riemannian methods which come with some theoretical guarantees are the Riemannian steepest-descent algorithm and the Riemannian conjugate gradient algorithm, both of them being known to converge to a critical point of Problem \ref{eq:BM} \cite{boumal2018global, sato2015new}. In practice however, many more algorithms for nonlinear optimization can be and have been generalized to manifold settings (though their convergence properties have not necessarily been studied) and are used to solve problems such as Problem \ref{eq:BM} (see, e.g., the list of solvers implemented in Manopt \cite{manopt}, one of the main toolboxes for optimization over manifolds). 

We conclude this section by mentioning a few applications where Riemannian optimization methods have recently been used. These include machine learning problems such as synchronization \cite{mei2017solving}, phase recovery \cite{sun2018geometric}, dictionary learning \cite{sun2016complete}, and low-rank matrix completion \cite{boumal2011rtrmc}; but also robotics problems such as the simultaneous localization and mapping (SLAM) problem (cf. Section \ref{sec:structure pop}). In \cite{carlone2015lagrangian, carlone2016planar, Rosen19}, the authors propose semidefinite programming relaxations for the maximum likelihood estimation (MLE) problem arising from SLAM. These relaxations provide exact solutions to the MLE problem assuming that noise in the measurements made by the robot are below a certain threshold (see \cite{Rosen19} for exact conditions). In \cite{Rosen19}, the authors tackle the challenge of scalability by demonstrating that the resulting SDP admits low-rank solutions and using a Burer-Monteiro factorization combined with Riemannian optimization techniques (e.g., a Riemannian trust-region approach) to solve the SDP. Thus, assuming that the conditions on the measurement noise are satisfied, one can find globally optimal solutions to the MLE problems arising from Pose-graph SLAM (and other related MLE problems arising from camera motion estimation and sensor network localization problems) via this method. 

\subsubsection{Coordinate descent approaches}

One of the more recent trends for solving problems of the form of Problem \ref{eq:BM} has been to use coordinate descent (or block-coordinate descent) methods \cite{erdogdu2018convergence,wang2017mixing}. They have been proposed to solve diagonally-constrained problems, i.e., problems of the type
\begin{equation}\label{eq:diag.constr}
\begin{aligned}
&\min_{V \in \mathbb{R}^{n \times r}} &&\mbox{Tr}(CVV^T)\\
&\text{s.t. } &&||v_i^T||=b_i, i=1,\ldots,n.
\end{aligned}
\end{equation}
where $||\cdot||$ is the 2-norm and $v_i^T$ is the $i^{th}$ row of $V$. This subset of SDPs cover machine learning applications such as certain formulations of graphical model inference and community detection \cite{erdogdu2018convergence}. Coordinate-descent methods are easy to explain for this case. Let us assume that $b_i=1$ for $i=1,\ldots,m$ for clarity of exposition. Initialization consists in randomly choosing $v_i, i=1,\ldots,n$ on the sphere. Then, at each iteration, an index $i$ is picked and fixed and the goal is to find a feasible vector $v_i$ such that $\mbox{Tr}(CVV^T)$ is minimized. 
We update the objective by replacing the previous $v_i$ by this new $v_i$ and proceed to the next iteration. This algorithm can be slightly modified to consider descent on blocks of coordinates, rather than one single coordinate at each time-step \cite{erdogdu2018convergence}. In that setting, it can be shown that if the indices are picked appropriately, then coordinate block-descent methods converge to an $\epsilon$-accurate critical point of Problem \ref{eq:diag.constr} in time $O(\frac{1}{\epsilon})$. Each iteration of coordinate descent methods is much less costly (by an order of $n$ approximately)  than its Riemannian trust-region counterpart however. As a consequence, numerical results tend to show that as $n$ increases, using coordinate descent methods can be preferable \cite{erdogdu2018convergence}.

\subsection{Frank-Wolfe based methods} \label{subsec:FW}

We start this subsection with a quick presentation of the Frank-Wolfe algorithm  \cite{frank1956algorithm}(also known as the conditional gradient algorithm), a version of which will be used in the approaches we present. This is an algorithm for constrained optimization problems of the type $\min_{x \in S} f(x)$, where $S$ is  a compact and convex set and $f$ is a convex and differentiable function. At iteration $k$, a linear approximation of $f$ (obtained using the first-order Taylor approximation of $f$ around the current iterate $x_k$) is minimized over $S$ and a minimizer $s_k$ is obtained. The next iterate $x_{k+1}$ is then constructed by taking a convex combination of the previous iterate $x_k$ and the minimizer $s_k$ before repeating the process. Intuitively, the reason why Frank-Wolfe algorithms are popular for obtaining low-rank feasible solutions for SDPs is due to the following fact: if the algorithm is initialized with a rank-1 matrix and if each minimizer is also a rank-1 matrix, then the $k^{th}$ iterate is of rank at most $k$. This enables us to control the rank of the solution that is given as output. We first present an algorithm developed by Hazan \cite{hazan2008sparse} in Section \ref{subsec:Hazan} and follow up with some extensions of Hazan's algorithm in Section \ref{subsec:other.FW}.

\subsubsection{Hazan's algorithm} \label{subsec:Hazan}
Hazan's algorithm is designed to produce low-rank solutions to SDPs of the type
\begin{equation}\label{eq:SDP.compact}
\begin{aligned}
&\min_{X \in S_n} &&\mbox{Tr}(CX)\\
&\text{s.t. } &&\mbox{Tr}(A_iX)=b_i,~i=1,\ldots,m\\
& &&X \succeq 0 \text{ and } \mbox{Tr}(X)=1.
\end{aligned}
\end{equation}
The additional constraint $\mbox{Tr}(X)=1$ is required by the algorithm, but can be slightly relaxed to the constraint $\mbox{Tr}(X)\leq t$, where $t$ is fixed \cite[Chapter 5]{gartner2012approximation}. We consider the version with the $\mbox{Tr}(X)=1$ constraint here. A key subroutine of Hazan's algorithm is solving the following optimization problem:
\begin{equation}\label{eq:Hazan}
\begin{aligned}
&\min_{X \in S_n} &&f(X)\\
&\text{s.t. } &&\mbox{Tr}(X)=1 \text{ and } X\succeq 0,
\end{aligned}
\end{equation}
where $f$ is a convex function. This is of interest in itself as some problems can be cast in the form \ref{eq:Hazan}, an example being given in \cite{hazan2008sparse}. Let $f^*$ be the optimal value of this problem. It is for this subroutine that a Frank-Wolfe type iterative algorithm is used. It can be described as follows: initalize $X_1$ to be a rank-$1$ matrix with trace one. Then, at iteration $k$, compute the eigenvector $v_k$ corresponding to the maximum eigenvalue of $\nabla f(X_k)$ written in matrix form. Finally, for a step size $\alpha_k$, set $X_{k+1}=(1-\alpha_k)X_k+\alpha_k (-v_kv_k^T)$ and iterate. This algorithm returns a $\frac{1}{k}$-approximate solution to \eqref{eq:Hazan} (i.e., a matrix $X$ such that $\mbox{Tr}(X) \leq 1+\frac{1}{k}$ and $f(X) \leq f^*+\Omega(\frac{1}{k})$) with rank at most $k$ in $k$ iterations. 

In order to see how this subroutine can be used to solve Problem \ref{eq:SDP.compact}, first note that one can reduce Problem \ref{eq:SDP.compact} to a sequence of SDP feasibility problems via the use of bisection. It is enough as a consequence to explain how one can leverage the Frank-Wolfe type algorithm described above to solve an SDP feasibility problem of the type: $$\mbox{Tr}(A_iX) \leq b_i,i=1,\ldots,m, X\succeq 0, \mbox{Tr}(X)=1.$$ 
This can be done by letting $f(X)=\frac{1}{M} \log\Big{(}\sum_{i=1}^m e^{M\cdot (\mbox{Tr}(A_iX)-b_i)}\Big{)}$ where $M=k\log(m)$, $\frac{1}{k}$ being the desired accuracy; see \cite{hazan2008sparse}. The algorithm is then guaranteed to return a $\frac{1}{k}$-approximate solution of rank at most $k$ in $O(k^2)$ iterations.

It is interesting to note that this latter algorithm is nearly identical to the multiplicative weights update algorithm for SDPs which also produces low-rank solutions \cite{arora2005fast} and has exactly the same guarantees. However it is derived in a completely different fashion.

\subsubsection{Other methods based on Frank-Wolfe} \label{subsec:other.FW}
A caveat of Hazan's algorithm is that the rank of the solution returned is linked to its accuracy. One may think that, if the solution is known to be low-rank and the Frank-Wolfe algorithm is initialized with a low-rank matrix, then all iterates of the algorithm produce a low-rank matrix. Unfortunately this is not always the case. In fact, what can be observed in practice is that the Frank-Wolfe algorithm typically returns iterates with increasing rank up to a certain point, where the rank decreases again, with the final solution returned being low-rank; see, e.g., the numerical experiments in \cite{yurtsever2017sketchy}. As a consequence, if one requires an accurate solution, one will likely have to deal with high-rank intermediate iterates, with the storage and computational problems that this entails.

In this section, we present two different methods designed to avoid this issue for optimization problems of the type:
\begin{equation}\label{eq:app.FW}
\begin{aligned}
&\min_{X \in \mathbb{R}^{n \times d}} f(\mathcal{A}X)\\
&\text{s.t. } ||X||_{*}\leq 1,
\end{aligned}
\end{equation}
where $f:\mathbb{R}^{m} \rightarrow \mathbb{R}$ is a convex and differentiable function, $||\cdot||_*$ denotes the nuclear norm, and $\mathcal{A}:\mathbb{R}^{n\times d} \rightarrow \mathbb{R}^m$ is a linear operator that maps $X$ to $(\mbox{Tr}(A_1X),\ldots,\mbox{Tr}(A_mX))^T$. The relaxation of the low-rank matrix completion problem given in \ref{eq:matrix.completion.relaxation} e.g. fits into this format. Similarly to Problem \ref{eq:Hazan}, the Frank-Wolfe algorithm can be applied to Problem \ref{eq:app.FW} quite readily. Initialization is done as before by setting $X_0$ to be some rank-1 matrix. At iteration $k$, one computes an update direction by finding a pair $(u_k,v_k)$ of singular vectors corresponding to the maximum singular value of $
\mathcal{A}^*( \nabla f(\mathcal{A}X_k))$ where $\mathcal{A}^*:\mathbb{R}^m \mapsto \mathbb{R}^{n \times d}$ is the adjoint operator to $\mathcal{A}$ given by $\mathcal{A}^*z=\sum_{i=1}^m z_i A_i$. We then construct a new iterate $X_{k+1}=(1-\alpha)X_k-\alpha u_kv_k^T$ and repeat.

The first method, described in \cite{freund2017extended}, is a modification of the Frank-Wolfe algorithm as described in the previous paragraph. The key idea is that at every iteration, the algorithm chooses between two types of steps. One type is the step described above, which is derived from the singular vectors of the matrix $\mathcal{A}^*( \nabla f(\mathcal{A}X_k))$, (a ``regular step" in the paper); the other is what is called an ``in-face step''. One way to take such a step is to minimize $\text{Tr}((\mathcal{A}^* \nabla f(AX_k))^TX)$ over $X$ in the minimal face of the nuclear norm unit ball to which $X_k$ belongs. The matrix $Z_k$ obtained from this minimization is then used to produce the next iterate: $X_{k+1}=(1-\alpha)X_k+\alpha Z_k$. We refer the reader to \cite{freund2017extended} to see how this can be done efficiently. The advantage of such a step is that the matrix $X_{k+1}$ constructed has rank no larger than $X_k$. If the choice between a regular step and an in-face step is appropriately made, it can still be the case that one has a $\frac{1}{k}$-approximate solution after $k$ iterations, though the rank of the $k^{th}$ iterate could be much less than $k$ depending on the number of ``in face'' steps used.

The second method given in \cite{yurtsever2017sketchy} describes an alternative way of dealing with possibly high rank of the intermediate iterates $X_k$ produced by the Frank-Wolfe algorithm applied to \eqref{eq:app.FW}. To get around this problem, the authors devise a ``dual'' formulation of the Frank-Wolfe algorithm described above where the primal iterates $X_k \in \mathbb{R}^{n \times d}$ are replaced by ``dual" iterates $z_k \in \mathbb{R}^m$, which are of smaller size. There is a need however to keep track of the iterates $X_k$ in some sense as the ultimate goal is to recover a solution to \eqref{eq:app.FW} at the end of the algorithm. To do this, the authors suggest a procedure based on \emph{sketching}. The procedure is as follows: two random matrices $\Omega \in \mathbb{R}^{d \times j}$ and $\Psi \in \mathbb{R}^{l \times n}$ are generated such that $j$ and $l$ are of order of the rank $r$ of the solution to Problem \ref{eq:app.FW}, and two new matrices $Y$ and $W$ are obtained as $Y=X \Omega$ and $W=\Psi X$. Note that the matrices $(Y,W)$ typically have smaller rank than the iterates $X_k$ that they represent. Any update that is made in the algorithm is then reflected on $(Y,W)$ rather than on the matrix $X$. At the end of the algorithm, one can reconstruct a rank $r$ solution to the problem from the updated pair $(Y,W)$; see \cite{yurtsever2017sketchy} for details. Note that with this way of proceeding, one need never store or utilize high-rank intermediate iterates. 
\section{SCALABILITY VIA ADMM AND AUGMENTED LAGRANGIAN METHODS}
\label{sec:first order}

An attractive alternative to using interior-point methods for solving SDPs is to use \emph{first-order} methods. Compared to interior-point methods, first-order methods can scale to significantly larger problem sizes, while trading off the accuracy of the resulting output (more precisely, first-order methods can require more time to achieve similar accuracy to interior-point methods). First-order methods are thus particularly attractive for applications that demand scalability, but do not require extremely accurate feasible solutions (e.g., some machine learning applications highlighted below). Here, we describe two recent approaches that involve first-order methods: (i) an approach based on the Alternating Direction Method of Multipliers (ADMM) \cite{ODonoghue16}, and (ii) an approach based on the augmented Lagrangian method \cite{Zhao10, yang2015sdpnal}. 

\subsection{Solving SDPs using ADMM} 

In \cite{ODonoghue16}, the authors present a first-order method for solving general cone programs, with SDPs as a special case. Here, we first describe the basic idea behind ADMM and then describe how this is applied to SDPs. We refer the reader to \cite{Boyd11} for a thorough exposition to ADMM. 

Consider a convex optimization problem of the form:
\begin{equation}\label{eq:admm problem}
\begin{aligned}
&\min_{x, z} &&f(x) + g(z)\\
&\text{s.t. } &&x=z.
\end{aligned}
\end{equation}
Here, the convex functions $f$ and $g$ can be nonsmooth and take on infinite values (e.g., to encode additional constraints). 
ADMM is a method for solving problems of the form \ref{eq:admm problem} (more generally, ADMM can solve problems where $x$ and $z$ have an affine relationship). ADMM operates by iteratively updating solutions $x^k$, $z^k$, and a dual variable $\lambda^k$ (associated with the constraint $x=z$) via the following steps:
\begin{enumerate}
\item $x^{k+1} = \underset{x}{\textrm{argmin}} \ \Big{(} f(x) + \frac{\rho}{2}\|x-z^k-\lambda^k\|_2^2 \Big{)}$ 
\item $z^{k+1} = \underset{z}{\textrm{argmin}} \ \Big{(} g(z) + \frac{\rho}{2}\|x^{k+1}-z-\lambda^k\|_2^2 \Big{)}$ 
\item $\lambda^{k+1} = \lambda^k - x^{k+1} + z^{k+1}$, 
\end{enumerate}
where $\rho > 0$ is a parameter. The initializations $z^0$ and $\lambda^0$ are usually taken to be 0, but can be arbitrary. Under mild technical conditions (see \cite{Boyd11}), the steps above converge. In particular, the cost $f(x^k) + g(z^k)$ converges to the optimal value of Problem \ref{eq:admm problem}, the residual $x^k - z^k$ converges to zero, and $\lambda^k$ converges to an optimal dual variable. 

In \cite{ODonoghue16}, the general ADMM procedure above is applied to SDPs. A key idea is to apply ADMM to the \emph{homogeneous self-dual embedding} of the primal SDP (Problem \ref{eq:SDP}) and its dual (Problem \ref{eq:dual SDP}). The homogeneous self-dual embedding converts the primal-dual pair of optimization problems into a single convex \emph{feasibility} problem. This is done by embedding the Karush-Kuhn-Tucker (KKT) optimality conditions associated with the primal and dual SDPs into a single system of equations and semidefinite constraints. If the original primal-dual pair has a feasible solution, it can be recovered from a solution to the embedding. If on the other hand, the original primal-dual pair is not solvable, the embedding allows one to produce a certificate of infeasibility. The homogeneous self-dual embedding is very commonly used in interior-point methods \cite{sedumi, mosek}. We refer the reader to \cite{ODonoghue16} for more details on this technique.

To summarize, the approach presented in \cite{ODonoghue16} consists of the following three steps:
\begin{enumerate}
\item Form the homogeneous self-dual embedding associated with the primal-dual SDP pair;
\item Express this embedding in the form \ref{eq:admm problem} required by ADMM;
\item Apply the ADMM steps to this problem.
\end{enumerate}

In \cite{ODonoghue16}, the authors also present techniques for efficiently implementing Step 3. These include techniques for (i) efficiently performing the projections required to implement the ADMM steps, (ii) scaling/preconditioning the problem data in order to improve convergence, and (iii) choosing stopping criteria. The resulting approach achieves significant speedups as compared to interior-point methods for a number of large-scale conic optimization problems (with some loss in accuracy of the solution). As an example, \cite{ODonoghue16} considers the problem of \emph{robust principal component analysis} \cite{Candes11} in machine learning. This problem corresponds to recovering the principal components of a data matrix even when a fraction of its entries are arbitrarily corrupted by noise. More specifically, suppose one is given a (large) data matrix $M$ that is known to be decomposable as
$$M = L + S,$$ 
where $L$ is a low-rank matrix and $S$ is sparse (but nothing further is known, e.g., the locations or number of nonzero entries of $S$). Interestingly, this problem can be solved using semidefinite programming under relatively mild conditions (see \cite{Candes11}). For this problem, the ADMM approach presented in \cite{ODonoghue16} is able to tackle large-scale instances where standard interior-point solvers (SeDuMi \cite{sedumi} and SDPT3 \cite{SDPT3}) run out of memory. Moreover, the approach provides significant (order of magnitude) speedups on smaller instances, while resulting in only a small loss in accuracy (less than $3 \times 10^{-4}$ reconstruction error of the low-rank matrix $L$). 

\subsection{Augmented Lagrangian methods}\label{subsec:Aug.Lag.Method}

Next, we discuss two closely-related approaches for solving SDPs using augmented Lagrangian methods: SDPNAL \cite{Zhao10}, and SDPNAL+ \cite{yang2015sdpnal}. SDPNAL operates on the dual SDP \ref{eq:dual SDP} by defining an \emph{augmented Lagrangian} function $L_\sigma: \mathbb{R}^m \times S_n \rightarrow \mathbb{R}$ as:
\begin{equation}
\label{eq:augmented lagrangian}
L_\sigma(y, X) = -b^T y + \frac{1}{2\sigma}(\| \Pi(X - \sigma(C - \sum_{i=1}^m y_i A_i )) \|^2 + \|X\|^2),
\end{equation}
where $\sigma > 0$ is a penalty parameter, $\|\cdot\|$ corresponds to the Frobenius norm, and $\Pi(\cdot)$ is the metric projection operator onto the set $S_n^+$ of $n \times n$ symmetric psd matrices. More precisely, $\Pi(Y)$ is the (unique) optimal solution to the following convex problem:
\begin{equation}
\begin{aligned}
&\min_{Z \in S_n^+} && \frac{1}{2} \|Z - Y \|. \\
\end{aligned}
\end{equation}
The augmented Lagrangian function is continuously differentiable (since $\|\Pi(\cdot)\|^2$ is continuously differentiable) and may be viewed as the (usual) Lagrangian function associated with a ``regularized" version of the SDP \ref{eq:dual SDP} (see \cite{Rockafellar76, Boyd11} for a more thorough exposition). 

The augmented Lagrangian method iteratively updates solutions to the dual Problem~\ref{eq:dual SDP} and the primal Problem \ref{eq:SDP} using the following steps:
\begin{enumerate}
\item $y^{k+1} = \underset{y}{\textrm{argmin}} \ L_{\sigma_k} (y, X^k)$;
\item $X^{k+1} = \Pi(X^k - \sigma_k(C - \sum_{i=1}^m y_i A_i ))$;
\item $\sigma_{k+1} = \rho \sigma_k$ for $\rho \geq 1$.
\end{enumerate}
In order to implement this method, we thus require the ability to solve the ``inner" optimization problems in Steps 1 and 2 above. The objective functions in these inner optimization problems are convex and continuously differentiable, but not twice continuously differentiable (see \cite{Zhao10}). The SDPNAL method presented in \cite{Zhao10} works by applying a semismooth Newton method combined with the conjugate gradient method to approximately solve the inner problems. Conditions under which the resulting iterates converge to the optimal solution are established in \cite{Zhao10}. In \cite{yang2015sdpnal}, the authors present SDPNAL+, which builds on the approach presented in \cite{Zhao10}. However, for certain kinds of SDPs (\emph{degenerate} SDPs; see \cite{yang2015sdpnal} for a formal definition), SDPNAL can encounter numerical difficulties. In \cite{yang2015sdpnal}, the authors present SDPNAL+, which tackles this issue by directly working with the primal SDP and using a different technique to solve the inner optimization problems that arise from the augmented Lagrangian method (specifically, a ``majorized" semismooth Newton method). 

In \cite{yang2015sdpnal}, the authors compare the performance of SDPNAL+ with two other first order methods on various numerical examples. In particular, \cite{yang2015sdpnal} considers SDP relaxations for the machine learning problem of $k$-means clustering. Specifically, the goal here is to partition a given set of data points into $k$ clusters such that the total sum-of-squared Euclidean distances from each data point to its assigned cluster center is minimized. In \cite{Peng07}, the authors present semidefinite relaxations for finding approximate solutions to this NP-hard problem. The numerical experiments presented in \cite{yang2015sdpnal} demonstrate that SDPNAL+ is able to solve the resulting SDPs with higher accuracy as compared to the other first-order methods while also resulting in faster running times (e.g., order of magnitude gains in running time on large instances of the clustering problem).

\section{TRADING OFF SCALABILITY WITH CONSERVATISM}
\label{sec:scalability_conservatism}

In this section, we discuss approaches that afford gains in scalability, but are potentially \emph{conservative}. We use the word ``conservative" to refer to guaranteed feasible solutions that may be suboptimal (in contrast to approaches that may produce points that slightly violate some constraints). In particular, here we discuss approaches that provide the user with a tuning parameter for \emph{trading off} scalability with conservatism. Conceptually, we classify such approaches into two: (i) ``special-purpose" approaches that leverage domain knowledge associated with the target application, and (ii) ``general-purpose" approaches that apply broadly across application domains.

As an example of an approach that falls under the first category, we highlight recent work on \emph{neural network verification} using semidefinite programming. Over the last few years, a large body of work has explored the susceptibility of neural networks to adversarial inputs (or more generally, uncertainty in the inputs) \cite{Szegedy13, Papernot16, Zheng16improving, Moosavi17}. For example, in the context of image classification, neural networks can be made to switch their classification labels by adding very small perturbations (imperceptible to humans) to the input image. Motivated by the potentially serious consequences of such fragilities in safety-critical applications (e.g., autonomous cars), there has been a recent explosion of work on verifying the robustness of neural networks to perturbations to the input \cite{Tjeng17, Wong17, Liu19, Pulina12, Xiang18, wang2018efficient} (see \cite{Liu19} for a more comprehensive literature review). Recently, approaches based on semidefinite programming have been proposed to tackle this challenge \cite{Raghunathan18, Qin19, Fazlyab19}. The challenge of scalability for this application is particularly pronounced since the number of decision variables in the resulting SDPs grows with the number of neurons in the network. In \cite{Fazlyab19}, the authors propose an approach for verifying neural networks by capturing their input-output relationships using quadratic constraints. This allows one to \emph{certify} that for a given set of inputs (e.g., a set of images obtained by making small perturbations to a training image), the output is contained within a specified set (e.g., outputs that assign the same label). The authors consider neural networks with different types of activation functions (e.g., ReLU, sigmoid, etc.) and propose different sets of quadratic constraints for each one. Importantly, by including or excluding different kinds of quadratic constraints, the approach allows one to \emph{trade off} scalability with conservatism. Moreover, by exploiting the modular structure of neural networks, the approach scales to neural networks with roughly five thousand neurons. The scalability afforded by the modular approach, however, comes at the cost of conservatism.

Next, we highlight approaches that fall under the second category mentioned above, i.e., ``general-purpose" approaches for trading off scalability with conservatism. 

\subsection{DSOS and SDSOS optimization}\label{subsec:dsos}

%
%
%

In~\cite{dsos_siaga}, the authors introduce ``DSOS and SDSOS optimization'' as linear programming and second-order cone programming-based alternatives to sum of squares and semidefinite optimization that allow for a trade off between computation time and solution quality. 
The following definitions are central to their framework.

\begin{definition} [dd and sdd matrices]
	A symmetric matrix $A=(a_{ij})$ is \emph{diagonally dominant} (dd) if $$a_{ii} \geq \sum_{j \neq i} |a_{ij}|$$ for all $i$. A symmetric matrix $A$ is \emph{scaled diagonally dominant} (sdd) if there exists a diagonal matrix $D$, with positive diagonal entries, such that $DAD$ is dd. 
\end{definition}

\begin{definition} [dsos and sdsos polynomials]
A polynomial $p\mathrel{\mathop:}=p(x)$ is said to be \emph{diagonally-dominant-sum-of-squares} (dsos) if it can be written as $p(x)=\sum_i \alpha_i q_i^2(x)$, for some scalars $\alpha_i\geq 0$ and some polynomials $q_i(x)$ that each involve at most two monomials with a coefficient of $\pm 1$. We say that a polynomial $p$ is \emph{scaled-diagonally-dominant-sum-of-squares} (sdsos) if it can be written as $p(x)=\sum_i q_i^2(x),$ for some polynomials $q_i$ that involve at most two monomials with an arbitrary coefficient. 
\end{definition}

\begin{figure}
\begin{center}
    \mbox{
      \subfigure[The set of parameters $a$ and $b$ for which the polynomial      
$p(x_1,x_2)=x_1^4+x_2^4+ax_1^3x_2+(1-\frac12 a-\frac12 b)x_1^2x_2^2+2bx_1x_2^3$ is dsos/sdsos/sos.]
      {\label{subfig:dsos.sdsos.sos}\scalebox{0.41}{\includegraphics{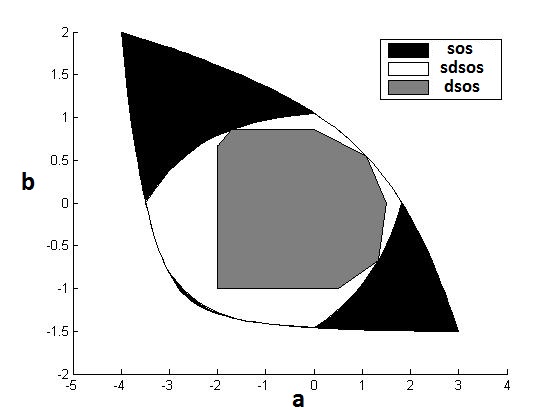}}} }
\mbox{\hspace{10mm}
      \subfigure[Reference~\cite{dsos_cdc14} uses SDSOS optimization to balance a model of a humanoid robot with 30 states and 14 control inputs on one foot.]
      {\label{subfig:atlas}\scalebox{0.43}{\includegraphics{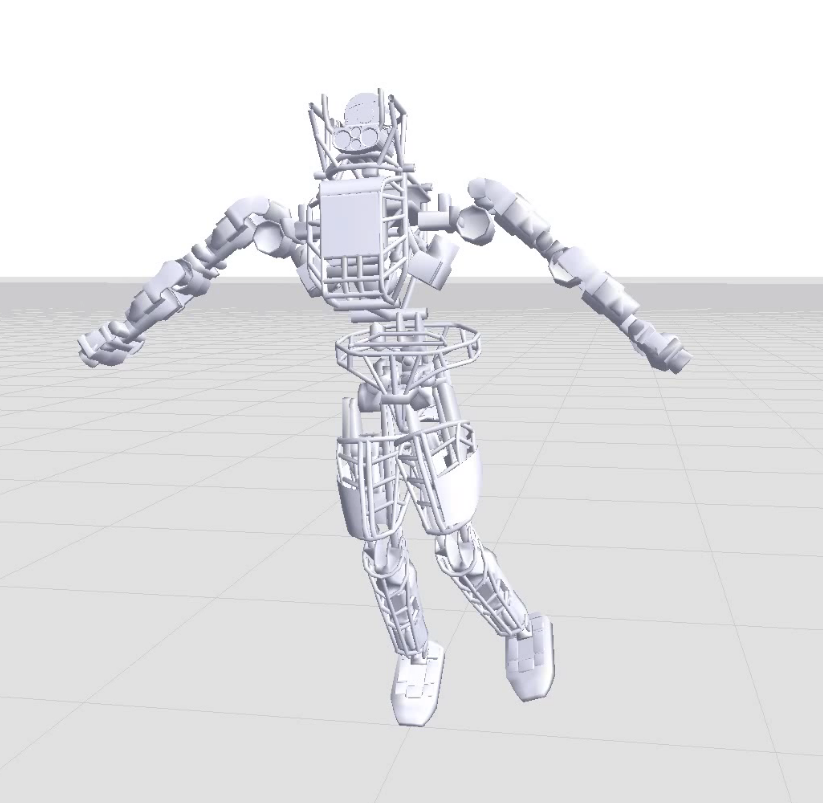}}}
} 
\vspace{10pt}
    \caption{Dsos and sdsos polynomials form structured subsets of sos polynomials which can prove useful when optimization over sos polynomials is too expensive.    
    }
\label{fig:dsos.and.atlas}
\end{center}
\vspace{-20pt}
\end{figure}

The following implications readily follow: $\mbox{dsos}\Rightarrow\mbox{sdsos}\Rightarrow\mbox{sos}.$ See Figure~\ref{subfig:dsos.sdsos.sos} for a comparison of the three notions on a parametric family of bivariate quartic polynomials. Similarly, in view of Gershgorin's circle theorem~\cite{gersh}, the implications $\mbox{dd}\Rightarrow\mbox{sdd}\Rightarrow\mbox{psd}$ are straightforward to establish. In~\cite{dsos_siaga}, the authors connect the above definitions via the following statement.

\begin{theorem}
A polynomial $p$ of degree $2d$ is dsos (resp. sdsos) if and only if it admits a representation as $p(x) = z^T(x)Qz(x)$, where $z(x)$ is the standard monomial vector of degree $\leq d$ and $Q$ is a dd (resp. sdd) matrix.
\end{theorem}

By combining this theorem with some linear algebraic observations, the authors show in~\cite{dsos_siaga} that optimization of a linear objective function over the intersection of the cone of dsos (resp. sdsos) polynomials of a given degree with an affine subspace can be carried out via linear programming (resp. second-order cone programming). These are two mature classes of convex optimization problems that can be solved significantly faster than semidefinite programs. The linear and second-order cone programs that arise from dsos/sdsos constraints on polynomials are termed \emph{``DSOS and SDSOS optimization problems''}. They are used to produce feasible, but possibly suboptimal, solutions to sum of squares optimization problems quickly. In the special case where the polynomials involved have degree two, DSOS and SDSOS optimization problems can be used to produce approximate solutions to semidefinite programs. We also remark that in applications where sum of squares programming is used as a \emph{relaxation}, i.e. to provide an outer approximation of a (typically intractable) feasible set, this approach replaces the semidefinite constraint with a membership constraint in the \emph{dual} of the cone of dsos or sdsos polynomials. See, e.g.,~\cite[Sect. 3.4]{isos_bp} for more details.

{\bf Impact on applications.} A key practical benefit of the DSOS/SDSOS approach is that it can be used as a ``plug-in'' in \emph{any} application area of SOS programming. The software package (cf. Section \ref{sec:software}) that accompanies reference~\cite{dsos_siaga} facilitates this procedure. In~\cite[Sect. 4]{dsos_siaga}, it is shown via numerical experiments from diverse application areas---polynomial optimization, statistics and machine learning, derivative pricing, and control theory---that with reasonable tradeoffs in accuracy, one can achieve noticeable speedups and handle problems at scales that are currently beyond the reach of traditional sum of squares approaches.\footnote{Comparisons are made in~\cite{dsos_siaga} with SDP solvers such as MOSEK, SeDuMi, and SDPNAL+.} For the problem of sparse principal component analysis in machine learning for example, the experiments in~\cite{dsos_siaga} show that the SDSOS approach is over 1000 times faster than the standard SDP approach on $100\times 100$ input instances while sacrificing only about 2 to 3\% in optimality; see~\cite[Sect. 4.5]{dsos_siaga}. As another example, Figure~\ref{subfig:atlas} illustrates a humanoid robot with 30 state variables and 14 control inputs that SDSOS optimization is able to stabilize on one foot~\cite{dsos_cdc14}. A nonlinear control problem of this scale is beyond the reach of standard solvers for sum of squares programs at the moment. In a different paper~\cite{dsos_opt_letters}, the authors show the potential of DSOS and SDSOS optimization for real-time applications. More specifically, they use these techniques to compute, every few milliseconds, certificates of collision avoidance for a simple model of an unmanned vehicle that navigates through a cluttered environment. 



{\bf Some guarantees of the DSOS/SDSOS approach.} On the theoretical front, DSOS and SDSOS optimization enjoy some of the same guarantees as those enjoyed by SOS optimization. For example, classical theorems in algebraic geometry can be utilized to conclude that any even positive definite homogeneous polynomial is the ratio of two dsos polynomials~\cite[Sect. 3.2]{dsos_siaga}. From this observation alone, one can design a hierarchy of linear programming relaxations that can solve any copositive program~\cite{dur2010copositive} to arbitrary accuracy~\cite[Sect. 4.2]{dsos_siaga}. This idea can be extended to achieve the same result for any polynomial optimization problem with a compact feasible set; see~\cite[Sect. 4.2]{POP_hierarchy}.

\subsection{Adaptive improvements to DSOS and SDSOS optimization}\label{subsec:dsos.improvements}

While DSOS and SDSOS techniques result in significant gains in terms of solving times and scalability, they inevitably lead to some loss in solution accuracy when compared to the SOS approach. In this subsection, we briefly outline two possible strategies to mitigate this loss. These strategies solve a sequence of adaptive linear programs (LPs) or second-order cone programs (SOCPs) that inner approximate the feasible set of a sum of squares program in a direction of interest.
%
For brevity of exposition, we explain how the strategies can be applied to approximate the generic semidefinite program given in Problem \ref{eq:SDP}.
A treatment tailored to the case of sum of squares programs can be found in the references we provide.


\subsubsection{Iterative change of basis}\label{subsubsec:cholesky} In~\cite{isos_bp}, the authors build on the notions of diagonal and scaled diagonal dominance to provide a sequence of improving inner approximations to the cone $P_n$ of psd matrices in the direction of the objective function of an SDP at hand. The idea is simple: define a family of cones\footnote{One can think of $DD(U)$ as the set of matrices that are dd after a change of coordinates via the matrix $U$.} $$DD(U)\mathrel{\mathop{:}}=\{M \in S_n~|~ M=U^TQU \text{ for some dd matrix } Q \},$$parametrized by an $n \times n$ matrix $U$. Optimizing over the set $DD(U)$ is an LP since $U$ is fixed, and the defining constraints are linear in the coefficients of the two unknown matrices $M$ and $Q$. Furthermore, the matrices in $DD(U)$ are all psd; i.e., $\forall U,$ $DD(U) \subseteq P_n$.

The proposal in~\cite{isos_bp} is to solve a sequence of LPs, indexed by $k$, by replacing the condition $X \succeq 0$ by $X \in DD(U_k)$:

\begin{flalign} 
\label{eq:LPChol}
			DSOS_k\mathrel{\mathop{:}}=\underset{X\in S_n}{\text{min}} \hspace*{1cm} & \mbox{Tr}(CX)  \\
			\text{s.t.} \hspace*{1cm} & \mbox{Tr}(A_iX)=b_i,\quad  i=1,\ldots,m, \nonumber \\
			& X\in DD(U_k). \nonumber
\end{flalign}


The sequence of matrices $\{U_k\}$ is defined as follows
\begin{equation}\label{eq:defUk}
\begin{aligned}
U_0&=I\\
U_{k+1}&=\text{chol}(X_k),
\end{aligned}
\end{equation}
where $X_k$ is an optimal solution to the LP in \ref{eq:LPChol} and chol(.) denotes the Cholesky decomposition of a matrix (this can also be replaced with the matrix square root operation).

Note that the first LP in the sequence optimizes over the set of diagonally dominant matrices. By defining $U_{k+1}$ as a Cholesky factor of $X_k$, improvement of the optimal value is guaranteed in each iteration. Indeed, as $X_k=U_{k+1}^T I U_{k+1}$, and the identity matrix $I$ is diagonally dominant, we see that $X_{k} \in DD(U_{k+1})$ and hence is feasible for iteration $k+1$. This implies that the optimal value at iteration $k+1$ is at least as good as the optimal value at the previous iteration; i.e., $DSOS_{k+1}\leq DSOS_k$ (in fact, the inequality is strict under mild assumptions; see~\cite{isos_bp}).

In an analogous fashion, one can construct a sequence of SOCPs that inner approximate $P_n$ with increasing quality. This time, the authors define a family of cones $${ SDD(U)\mathrel{\mathop{:}}=\{M \in S_n ~|~ M=U^TQU, \text{ for some sdd matrix } Q\},}$$ parameterized again by an $n\times n$ matrix $U$. For any $U$, optimizing over the set $SDD(U)$ is an SOCP and we have $SDD(U)\subseteq P_n$. This leads us to the following iterative SOCP sequence:

\begin{flalign} 
\label{eq:SOCPchol}
			SDSOS_k\mathrel{\mathop{:}}=\underset{X\in S_n}{\text{min}} \hspace*{1cm} & \mbox{Tr}(CX)  \\
			\text{s.t.} \hspace*{1cm} & \mbox{Tr}(A_iX)=b_i,\quad  i=1,\ldots,m, \nonumber \\
			& X\in SDD(U_k). \nonumber
\end{flalign}


Assuming existence of an optimal solution $X_k$ at each iteration, we can define the sequence $\{U_k\}$ iteratively in the same way as was done in Equation \ref{eq:defUk}. Using similar reasoning, we have $SDSOS_{k+1}\leq SDSOS_k$. In practice, the sequence of upper bounds $\{SDSOS_k\}$ approaches faster to the SDP optimal value than the sequence of the LP upper bounds $\{DSOS_k\}$.
%

\begin{figure}[h!]
	\begin{center}
		\mbox{
			\subfigure[LP inner approximations]
			{\label{subfig:DDCones}\scalebox{0.4}{\includegraphics{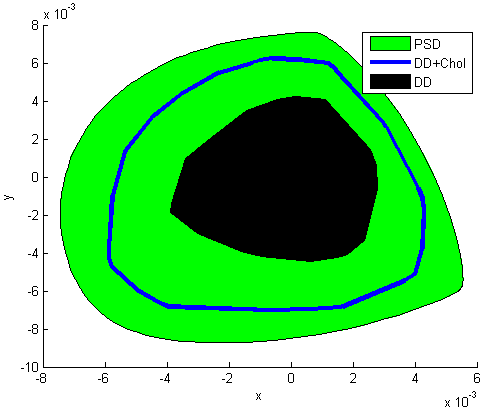}}}}
		\mbox{
			\subfigure[SOCP inner approximations]
			{\label{subfig:SDDCones}\scalebox{0.4}{\includegraphics{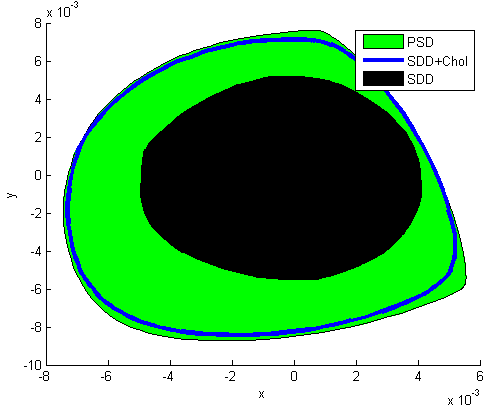}}}
		}
		
		\vspace{15pt}
		\caption{Figure reproduced from~\cite{isos_bp} showing improvement (in all directions) after one iteration of the change of basis algorithm.}
		\label{fig:DDSDDCones}
	\end{center}
	\vspace{10pt}
\end{figure}

Figure~\ref{fig:DDSDDCones} shows the improvement (in every direction) obtained just by a single iteration of this approach. The outer set in green in both subfigures is the feasible set of a randomly generated semidefinite program. The sets in black are the diagonally dominant (left) and the scaled diagonally dominant (right) inner approximations. What is shown in dark blue in both cases is the boundary of the improved inner approximation after one iteration. Note that the SOCP in particular fills up almost the entire spectrahedron in a single iteration.


\subsubsection{Column generation}\label{subsubsec:column.generation}
In~\cite{isos_cg}, the authors design another iterative method for inner approximating the set of psd matrices using linear and second order cone programming. Their approach combines DSOS/SDSOS techniques with ideas from the theory of column generation in large-scale linear and integer programming. The high-level idea is to approximate the SDP in Problem \ref{eq:SDP} by a sequence of LPs (parameterized by $t$):
\begin{flalign} 
\label{eq:LP.colgen}
			\underset{X\in {S}_n, \alpha_i}{\text{min}} \hspace*{1cm} & \mbox{Tr}(CX)  \\
			\text{s.t.} \hspace*{1cm} & \mbox{Tr}(A_iX)=b_i,\quad  i=1,\ldots,m, \nonumber \\
			& X=\sum_{i=1}^t \alpha_i B_i, \nonumber \\
			\ & \alpha_i\geq 0, \quad i=1,\ldots, t, \nonumber
\end{flalign}
where $B_1,\ldots,B_t$ are fixed psd matrices. These matrices are initialized to be the extreme rays of $n\times n$ dd matrices which turn out to be all rank one matrices $v_iv_i^T$, where the vector $v_i$ has at most two nonzero components, each equal to $\pm 1$~\cite{dd_extreme_rays}. Once this initial LP is solved, one adds one (or sometimes several) new psd matrices $B_j$ to Problem \ref{eq:LP.colgen} and resolves. This process then continues. In each step, the new matrices $B_j$ are picked carefully to bring the optimal value of the LP closer to that of the SDP. Usually, the construction of $B_j$ involves solving a ``\emph{pricing subproblem}'' (in the terminology of the column generation literature), which adds appropriate cutting planes to the dual of Problem \ref{eq:LP.colgen}; see~\cite{isos_cg} for more details.

\begin{figure}[t]
	\begin{center}
		\mbox{
			\subfigure[LP iterations]
			{\label{subfig:dsos.iters}\scalebox{0.4}{\includegraphics{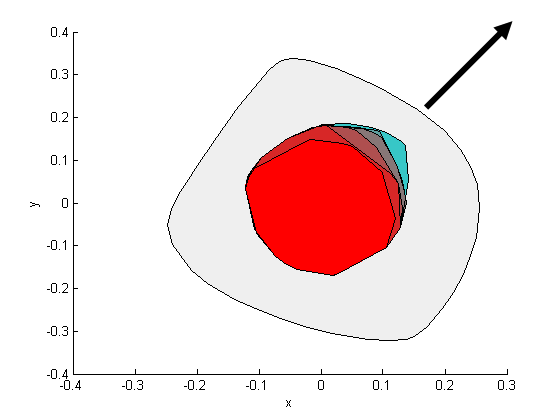}}}}
		\mbox{
			\subfigure[SOCP iterations]
			{\label{subfig:sdsos.iters}\scalebox{0.4}{\includegraphics{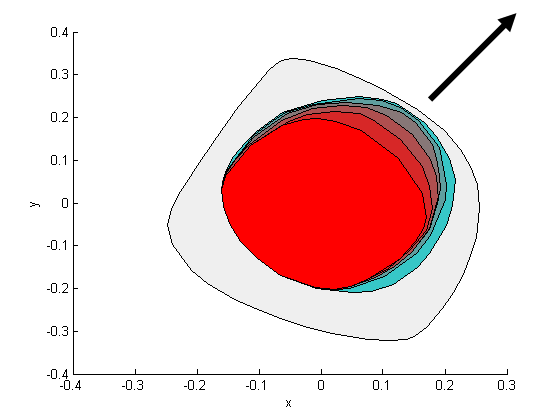}}}
		}
		\vspace{15pt}
		\caption{Figure reproduced from~\cite{isos_cg} showing the successive improvement on the dd (left) and sdd (right) inner approximation of the feasible set of a random SDP via five iterations of the column generation method.}
		\label{fig:colgen}
	\end{center}
\end{figure}

The SOCP analog of this process is similar. The SDP in Problem \ref{eq:SDP} is inner approximated by a sequence of SOCPs (parameterized by $t$):
\begin{flalign} 
\label{eq:SOCP.colgen}
			\underset{X\in{S}_n, \Lambda_i\in{S}_2}{\text{min}} \hspace*{1cm} & \mbox{Tr}(CX)  \\
			\text{s.t.} \hspace*{1cm} & \mbox{Tr}(A_iX)=b_i,\quad  i=1,\ldots,m, \nonumber \\
			& X=\sum_{i=1}^t V_i\Lambda_i V_i^T, \nonumber  \\
						\ & \Lambda_i\succeq 0, \quad i=1,\ldots, t, \nonumber 
\end{flalign}
where $V_1,\ldots,V_t$ are fixed $n \times 2$ matrices. They are initialized as the set of matrices that have zeros everywhere, except for a 1 in the first column in position $j$ and a 1 in the second column in position $k\neq j$. This gives exactly the set of $n\times n$ sdd matrices~\cite[Lemma 3.8]{dsos_siaga}. In subsequent steps, one (or sometimes several) appropriately-chosen matrices $V_i$ are added to Problem \ref{eq:SOCP.colgen}. These matrices are obtained by solving pricing subproblems and help bring the optimal value of the SOCP closer to that of the SDP in each iteration.

%

Figure~\ref{fig:colgen} shows the improvement obtained by five iterations of this process on a randomly generated SDP, where the objective is to maximize a linear function in the northeast direction.

\newpage
\section{SOFTWARE: SOLVERS AND MODELING LANGUAGES}
\label{sec:software}

In this section, we list software tools that allow a user to specify and solve SDPs. This list is not meant to be exhaustive; rather, the goal is to provide the interested practitioner with a starting point in terms of software for implementing some of the approaches reviewed in this paper. 

\begin{issues}[SOFTWARE TOOLS]
\begin{itemize}
\item Software for specifying SDPs: \\
YALMIP \cite{yalmip}, CVXPY \cite{cvxpy}, PICOS \cite{picos}, SPOT \cite{spot}, JuMP \cite{jump}.   \\
These software packages allow one to easily specify SDPs and provide interfaces to various solvers listed above. YALMIP and SPOT are MATLAB packages, CVXPY and PICOS are Python packages, and JuMP is a Julia package.  
\item General-purpose SDP solvers that implement interior-point methods: \\
MOSEK \cite{mosek}, SeDuMi \cite{sedumi}, SDPT3 \cite{SDPT3}. \\
The MOSEK solver provides an interface with MATLAB, C, Python, and Java. SeDuMi and SDPT3 interface with MATLAB. 
\item Solvers based on ADMM and augmented Lagrangian methods (cf. Section \ref{sec:first order}): \\
SCS \cite{scs}, SDPNAL/SDPNAL+ \cite{yang2015sdpnal}. \\
SDPNAL/SDPNAL+ provide a MATLAB interface, while SCS provides interfaces to C, C++, Python, MATLAB, R, and Julia.  
\item Software for Riemannian Optimization (cf. Section \ref{subsec:BM}): \\
Manopt \cite{manopt}, Pymanopt \cite{pymanopt}. \\
Manopt is a MATLAB-based toolbox for optimization on manifolds, while Pymanopt is Python-based.
\item Software for exploiting chordal sparsity in SDPs (cf. Section \ref{sec:sparsity}): \\
SparseCoLO \cite{Fujisawa09}, CDCS \cite{Zheng19}. \\
SparseCoLO is a MATLAB toolbox for converting SDPs with chordal sparsity to equivalent SDPs with smaller-sized SDP constraints. CDCS is a solver that exploits chordal sparsity and provides a MATLAB interface. 
\item Software for exploiting structure via facial reduction (cf. Section \ref{sec:facial reduction}): \\
frlib \cite{Permenter18}. \\
This is a MATLAB toolbox and interfaces with the YALMIP and SPOT packages mentioned above. 
\item Parsers for SOS programs (cf. Section \ref{subsec:sos}): \\
YALMIP \cite{yalmip}, SPOT \cite{spot}, SOSTOOLS \cite{sostools}, SumofSquares.jl \cite{sumofsquares-julia}. \\
These packages allow one to specify SOS programs and convert them to a form that can be solved using SDP solvers. YALMIP, SPOT, and SOSTOOLS are MATLAB packages, while SumofSquares.jl is a Julia package. 
\item Software for parsing DSOS and SDSOS programs (cf. Section \ref{subsec:dsos}): \\
SPOT \cite{spot}. \\
This MATLAB toolbox allows one to parse DSOS and SDSOS programs. SPOT provides interfaces to LP and SOCP solvers including MOSEK. See the appendix of \cite{dsos_siaga} for a tutorial on solving DSOS and SDSOS programs using SPOT.
\end{itemize}
\end{issues}
\section{CONCLUSIONS}
\label{sec:conclusions}

Semidefinite programming is an exciting and active area of research with a vast number of applications in fields including machine learning, control theory, and robotics. Historically, the lack of scalability of semidefinite programming has been a major impediment to fulfilling the potential that it brings to these application domains. In this paper, we have reviewed recent approaches for addressing this challenge including (i) techniques that exploit structure such as sparsity and symmetry in a problem, (ii) approaches that produce low-rank feasible solutions to SDPs, (iii) methods for solving SDPs via augmented Lagrangian and ADMM techniques, and (iv) approaches that trade off scalability with conservatism, including techniques for approximating SDPs with LPs or SOCPs. We have also provided a list of software packages associated with these approaches. Our hope is that this paper will serve as an entry-point for both practitioners working in application domains that demand scalability, and researchers seeking to contribute to theory and algorithms for scalable semidefinite programming. 

Exciting and significant work remains to be done on building upon the advancements reviewed here. Semidefinite programming is still far from being a mature technology similar to linear or quadratic programming. Potential directions for future work include: (i) better theoretical understanding of the convergence properties of algorithms for low-rank SDPs (cf. Section \ref{sec:low.rank}),  (ii) identifying structure beyond chordal sparsity, symmetry, and degeneracy (cf. Section \ref{sec:structure}) that arise in practice and methods for exploiting such structure, (iii) exploring different first-order methods for solving SDPs and understanding their relative merits (cf. Section \ref{sec:first order}), (iv) understanding the power of LPs and SOCPs for approximating SDPs in an adaptive (as opposed to one-shot) fashion (cf. Section \ref{sec:scalability_conservatism}), (v) finding ways to combine the different approaches for improving scalability reviewed in this paper, and (vi) further developing mature software that allows practitioners to deploy semidefinite programming technology on future applications including ones that involve solving SDPs in real-time.

%
%

\section*{DISCLOSURE STATEMENT}
The authors are not aware of any affiliations, memberships, funding, or financial holdings that
might be perceived as affecting the objectivity of this review. 

\section*{ACKNOWLEDGMENTS}
We thank Cemil Dibek for his constructive feedback on the first draft of this manuscript. This work is partially supported by the DARPA Young Faculty Award, the CAREER Award of the NSF, the Google Faculty Award, the Innovation Award of the School of Engineering and Applied Sciences at Princeton University, the Sloan Fellowship, the MURI Award of the AFOSR, the Office of Naval Research [Award Number: N00014-18-1-2873], the NSF CRII award [IIS-1755038], and the Amazon Research Award. 


\bibliographystyle{abbrvnat}
\bibliography{bib-sdp-scalability} 




\end{document}